\documentclass[a4paper,12pt]{amsart}
\usepackage{amssymb}
\usepackage{ifthen}
\usepackage[dvips]{graphicx}
\newtheorem{lem}{Lemma}
\newtheorem{thm}{Theorem}
\newtheorem{rk}{Remark}

\begin{document}
\bibliographystyle{amsplain}
\title{ a  Problem of Bombieri on\\ Univalent
 Functions}
\author{D. Aharonov  and D. Bshouty}
\begin{abstract}The famous Bieberbach Conjecture from 1916 on the coefficients of normalized univalent functions defined in the unit disk \cite{bil} that was finally proved by de Branges \cite{dbl} some 70 years later, drifted many complex analysts attention to other subjects. Those who continued to explore de Branges method and push it as far as possible were not aware of where it may lead. Surprisingly enough, a paper that fell in our hands \cite{dxh} contained a way to tackle one of the problems of Bombieri \cite{boe} on the  behavior of the coefficients of univalent functions. We shall give an account of the history of the problem and a revised version of it.
\end{abstract}
\maketitle


\vspace{.3cm}

\section {Introduction}

\thanks{}

\pagestyle{myheadings}

\markboth{}{}

Since 1916, The Bieberbach Conjecture (BC) \cite{bil} , was the basic open problem of geometric function theory till its proof by de Branges in 1984 \cite{dbl}.
This conjecture has a long history. It states that for functions in the class $\mathbf{S}$ of one-to-one analytic functions $f$ defined on the unit disk $D$ and normalized by
$$f(z)=z+\sum_{n=2}^\infty a_nz^n;\;\;z\in D $$
the relation
$$|a_n|\leq n;\;\; n\geq 2$$
holds true with equality only for the Koebe function
$$ f(z)=\frac z{(1-z)^2}=z+\sum_{n=2}^\infty nz^n,$$
or one of its rotations. Bieberbach \cite{bil} himself proved it for $n=2$ by appealing to the area theorem for the class $\Sigma$ of normalized univalent functions in the exterior of  $D,\; D^c.$ Subsequently, in  1923 L\"{o}wner \cite{lok} introduced his partial differential chain approach and gave the first proof of the conjecture for $n=3.$
When the method did not seem to meet the expectations to solve the conjecture, various techniques were introduced. Powerful variational methods were introduced in 1938 by Schiffer and, Duren and Schiffer \cite{scm1,scm2,DS} which proved useful to attack several extremal problems in the class $\mathbf{S}.$ In particular,in 1955, with tedious work Garabedian and Schiffer \cite{GS1} proved  BC for $n=4$  and again in 1960 the proof was simplified by Charzy\'{n}ski and Schiffer \cite{CS}. Subsequently, other proofs were further simplified using other techniques.  In 1939 the Grunsky inequalities \cite{grh} were the first necessary and and sufficient condition on the coefficients of an analytic function $g,$ defined in $D^c,$ to belong to $\Sigma$. They were based on the simple observation that $g$ is univalent if, and only if,
\begin{equation}
\log\frac{g(z)-g(\zeta)}{z-\zeta},\;\; \mbox{is analytic in} \;\; |z|>1,\; |\zeta|>1
\end{equation}
and the area theorem. These inequalities were generalized by Garabedian and Schiffer  \cite{GS} and in conjunction with the variational techniques lead  to the prove of BC for $n=5$ and $6.$
The belief in the truth of BC drew the attention of researchers to the investigation of the conjecture "near" the supposed extremal mapping, the Koebe function. In 1965 and 1967 the local BC was proved via the aformentioned inequalities for even and then odd $n$ \cite{GRS,GS} in the form: For each $n$ there exist $\epsilon_n>0$ such that
$$ \Re\{n- a_n\}>0 \;\;\mbox{ whenever } \;\;\Re\{2-a_2\}<\epsilon \;\; \mbox {for}\;\; 0<\epsilon<\epsilon_n.  $$

\vspace{.3cm}
 Few years earlier, in 1963 Bombieri, in his thesis, ( see \cite{haw} problem 6.3) proved  that there exist positive constants $c_n$ such that for every
 $f\in \mathbf S$
$$\big|\Re\{n-a_n\}\big |<b_n\;\Re\{2-a_2\} $$
and questioned about the size of the constants $b_n$. Moreover, he asked  whether  it is true that there exist positive constants $d_n$ such that
\begin{equation}
\big|n-|a_n|\big |<d_n\big(2-|a_2|\big).
\end{equation}
In 1967, Bombieri \cite{boe}, using the L\"owner method and the variational method of Duren and Schiffer \cite{DS},  proved a striking form of the local BC. Namely, there exist positive constants $e_n$ and $g_n$ such that
\begin{equation}
 \varliminf_{a_2\rightarrow 2} \frac {n-\Re\{a_n\}}{2-\Re\{a_2\}}\geq e_n,\;n-\mbox{even}\;;\;\;
 \varliminf_{a_3\rightarrow 3} \frac {n-\Re\{a_n\}}{3-\Re\{a_3\}}\geq g_n,\;n-\mbox{odd}
\end{equation}
and conjectured that for $n>m>1$
\begin{equation}
 B_{m,n}=\varliminf_{a_m\rightarrow m} \frac {n-\Re\{a_n\}}{m-\Re\{a_m\}} =\min_\theta\frac{n-\frac{\sin(n\theta)}{\sin(\theta)}}{m-\frac{\sin(m\theta)}{\sin(\theta)}}.
 \end{equation}
 
An explanation is due here. At the time, proving  BC seemed unaccessible except for the first few coefficients,  thus it was obvious to draw the attention to the local problem near the extremal function and then advance our understanding on the behaviour of all the coefficients there. Thus , for example,(3) can be interpreted to saying that The collection of all the tuples
 $$V_n=\{(\Re\{a_2\},\Re\{a_n\}),\; f(z)\in \mathbf S,\; n\; \mbox{even}\} $$
  near the point$(2,n)$ is bounded from above by a straight line through that point with a positive slope. As opposed to this, Bombieri's thesis result is  weaker  in that it bounds $V_n$ from above equally but with a negative slope as well as of a positive slope from below, but for all values of $0\leq|a_2|\leq2.$ The obvious Question for
$$ W_n=\{(|a_2|,|a_n|),\; f(z)\in \mathbf S,\; n\; \mbox{even}\} $$ 
arose. All these results being qualitative, quantitative results were the next step.

As to inequality (2), it follows from (3) \cite{had}. A direct proof based on L\"owner chains solely is found in \cite{bsd2}. It is based on  L\"owner's formula for $a_n(t)$ as represented in terms of the earlier coefficients and an inductive procedure.

Turning to quantitative results, inequality (2) in the subclass of starlike functions in  $\mathbf{S}$ is found in \cite{huj} with the  quantitative bound $d_n= n(n^2-1)/6$ which is asymptotically exact. Inequality (4) holds true for the subclass of $\mathbf{S}$ with real coefficients and for analytic variations of the Koebe mapping in $\mathbf S$ \cite{BH,PR}, however, it is not true for general variations. One example is  $B_{2,3}$ \cite{grr} where the exact value is also found, and others are  $ B_{2,4}$ and $B_{3,4}$ \cite{PV}.  If we were to summarize the foregoing, it boils down to the following: The L\"owner method and the Grunsky inequalities and its generalization were the most powerful tools in tackling the local BC  \cite{fic2}, and more so the Bombieri inequalities. In the  whole class $\mathbf{S}$ we lack of any quantitative information of that sort.

Tackling BC seemed to need a closer look either from the point of view of (1)  or else the L\"owner method.  Lebedev and Milin  developed a new method based on the Grunsky approach. Their idea was  to  estimate  the coefficients of a univalent function from those of its logarithm  \cite{M} (see also \cite{P} p. 78). In 1967  Milin \cite{mii} applied it to the relation
\begin{equation} \frac{f(z)}{z}=\exp\big[\sum_1^\infty 2\gamma_nz^n\big]\nonumber
\end{equation}
for odd  $f\in \mathbf S$ to show that $|a_n|<1.243n.$ In 1971 Lebedev and Milin ( see \cite{M,A}) showed that BC would follow from the inequalities
\begin{equation}
M_{n}(f)=\frac1n\sum_{k=1}^{n-1}(n-k)\bigg ( k|\gamma_k|^2-  \frac{1}k\bigg ) \leq 0,\;\;n=2,3,\dots,
\end{equation}
and conjectured the truth of (5). Using Lebedev-Milin inequalities, Aharonov  \cite{ahd} showed that $|a_n|<n$ whenever $|a_2|<1.05$
In 1972 FitzGerald \cite{fic1} applied another purely algebraic exponentiation technique to Goluzin inequalities
to prove that $|a_n|<1.081n.$ Subsequently, the  FitzGerald inequalities were used to improve Aharonov's result \cite{ehg1,bsd1}. The methods used for these results rest on using sharpened inequalities of the original work of Milin and FitzGerald . These ideas are efficient in delivering quantitative estimates for $d_n$ in (2) but only for small values of $|a_2|.$

The proof of de Branges \cite{dbl} of BC was directed to inequality (5) from which BC follows. Inequality (2) was thus taken care of from one side in the form $n-|a_n|\geq 0.$ At this stage the question that arose was to find a better linear or otherwise upper bound of the same sort. In view of (3), we have evidence that even and odd coefficients behave differently near $|a_2|=2$ and more so due to equation (1.4) in \cite{boe}, namely,

\begin{equation}
\varliminf_{a_2\rightarrow 2} \frac {3-\Re\{a_3\}}{(2-\Re\{a_2\})^{\frac32}}=\frac 83.\nonumber
\end{equation}

\vspace{.1cm}
\noindent
An improved inequality of (5) due to Dong \cite{dxh} is what we use to improve the upper bound of $|a_n|$ in terms of $|a_2|.$ We shall prove

\vspace{.4cm}
\noindent
{\bf Main Theorem.} {\it Let $f(z)\in \mathbf S. $  Then for  $m\in \mathbb N$ we have
$$\frac{2m-|a_{2m}|}{2-|a_2|}\geq\frac{3m}{8(4m^2-1)},$$
$$\frac{2m+1-|a_{2m+1}|}{(2-|a_2|)^{\frac32}}\geq \frac{3m(m+1)}{(2m-1)(2m+3)}.$$
 }

As to the other side of inequality (2),
namely,
\begin{equation}
\frac{n-|a_n|}{2-|a_2|}\leq d_n\nonumber
\end{equation}
it  amounts to looking for a linear lower bound of $|a_n|$ in terms of $|a_2|.$ In this respect,for functions with real coefficients in $\mathbf {S},$ the extremal mapping for $a_2=c; \; -2\leq c \leq 2$  is $f_c(z)=\frac z{1-2cz+z^2}$ \cite{BH}. In the class  $\mathbf {S}$ the sharp bounds for $|a_3|, |a_4|$ and $ |a_5|$ were found in specific intervals near $|a_2|=2$ where again  $f_c(z)$ is the extremal mapping  \cite{gor}. However to our knowledge finding  specific values of $d_n$ in $\mathbf {S}$ remains open.

\section {Preliminary results and proof of Main Theorem}

de Branges proof of Lebedev-Milin inequalities is based on a system of ordinary differential equations with specific initial conditions that produced nonincreasing solutions, now termed de Branges system. To improve  inequality (5), it was necessary to decrease the initial conditions while the solutions remain nonincreasing. This is presicely Dong's motivation and a practical example of this basic idea can be found in Li  (\cite{lij} p. 167).

\vspace{.3cm}
\begin{thm} Dong \cite{dxh}  Let $f(z)=z+\sum_{n=2}^\infty a_nz^n\in \mathbf S $ and set
\begin{equation}
 \frac12\log\frac{f(z)}{z}=\sum_{n=1}^\infty \gamma_nz^n.
 \end{equation}
Then for each $n=2,3,\dots $ we have
$$M_n(f)+K_n(f)\leq 0$$
where
\begin{equation}
K_{2m}(f)=\frac{3}{4(4m^2-1)}(1-|\gamma_1|^2),
\end{equation}
\begin{equation}
K_{2m+1}(f)=\frac{15m(m+1)}{(4m^2-1)(2m+3)}\big(\frac54-|\gamma_1|^2-|\gamma_2|^2\big).
\end{equation}
\end{thm}

\vspace{.3cm}
\begin{rk} Note that $ M_2(f)=-\frac12(1-|\gamma_1|^2)$ and that $ M_3(f)=-\frac23(\frac54-|\gamma_1|^2-|\gamma_2|^2\big).$
\end{rk}

\vspace{.3cm}
We shall need the following Lemmas
\begin{lem} (Lebedev-Milin Inequality \cite{M}) Let the formal power series $$g(z)= \sum_{k=1}^\infty \beta_kz^k$$ be given. Set
$$\exp\big(g(z)\big)=\sum_{k=0}^\infty p_kz^k.$$
 Then the sequence
\begin{equation}
Q_n=\frac{\exp[-M_n(g)]}n\sum_{k=0}^{n-1} |p_k|^2;\;\; n\geq 2
\end{equation}
is a monotone decreasing sequence.
\end{lem}

\noindent
\begin{lem} (Lemma 1, p.196 \cite {gol}) Let $\lambda(t);\; t\geq0 $ be an arbitrary continuous real function
except possibly for a finite number of discontinuities of the first kind. Suppose that $|\lambda(t)|\leq e^{-t}; \;t\geq0.$
Then by setting
$$\int_0^\infty\lambda^2(t)dt=(\nu+\frac12)e^{-2\nu};\; \nu\geq0\;\;\mbox{we have}\;\;\bigg| \int_0^\infty\lambda(t)dt\bigg|\leq (\nu+1)e^{-\nu}.$$
Equality holds only for functions $\lambda(t)=\pm \mu(t),$ where $\mu(t)=e^{-\nu}; \; 0<t<\nu$ and $\mu(t)=e^{-t};\; t>\nu.$
\end{lem}

\vspace{.3cm}
\noindent
\begin{rk} From Lemma 2 and the fact that the functions $\phi_1(\nu)=(\nu+1)e^{-\nu}$ and $\phi_2(\nu)=(\nu+\frac12)e^{-2\nu}$ are
monotonically decreasing we have the following: If
$$\bigg| \int_0^\infty\lambda(t)dt\bigg|= (\nu+1)e^{-\nu},\;\;\mbox{ then } \;\;
\int_0^\infty\lambda^2(t)dt\geq(\nu+\frac12)e^{-2\nu}.$$
Equality holds if, and only if, $\lambda(t)$ is as in the case of equality in Lemma 2.
\end{rk}

\vspace{.3cm}
\noindent
We recall the following relations
$$ \gamma_1=\frac{a_2}2 \;\; \mbox{and}\;\; \gamma_2=\frac{2a_3-a_2^2}4,$$ and  prove

\begin{lem} For $f(z)=\sum_{k=1}^\infty a_n z^n \in \mathbf S$ we have
\begin{equation}  |a_3-\frac {a_2^2}2|^2+|a_2|^2-5 =4(|\gamma_2|^2+|\gamma_1|^2-\frac54)\leq-\sqrt{2}(2-|a_2|)^{3/2}. \end{equation}
\end{lem}
\noindent
{\bf Proof: }
Let $f\in S, \;\; f(z)= z + a_2 z^2+a_3 z^3 +\dots;\;\;  z\in D.$ We  follow the proof  of Fekete and Szeg\"o's theorem on p.198 in \cite {gol}.
Using a certain rotation of $a_2$ we may assume  that $a_3-a_2^2\geq 0.$  For $\alpha=\frac12,$ equation (7) [{\it loc. cit.}],  yields
$$a_3-\frac12a_2^2=1-4\int_0^\infty e^{-2t}\lambda(t)^2dt+2\bigg(\int_0^\infty e^{-t}\lambda(t)dt\bigg)^2$$
where for some $\theta_0,\;\;\Re\{a_2e^{i\theta_0}\}=2\int_0^\infty e^{-t}\lambda(t)dt;\;\; \lambda(t)$ is a real function and $|\lambda(t)|\leq 1.$
By  Remark 2, there exist a $\nu\geq0$ such that for
$$|\Re\{a_2e^{i\theta_0}\}|=2 \bigg|\int_0^\infty e^{-t}\lambda(t)dt\bigg|= 2 (\nu+1)e^{-\nu};\;\nu>0$$
we have
$$\int_0^\infty e^{-2t}\lambda(t)^2dt\geq(\nu+\frac12 )e^{-2\nu}.$$
We conclude that
\begin{equation}
 |a_3-\frac12a_2^2|\leq 2 e^{-2 \nu}  \nu^2+1.\nonumber
 \end{equation}
\noindent
This inequality is invariant to rotations of  $a_2e^{i\theta_0}$ and therefore we conclude that
the same holds for  $|a_2|=2(\nu+1)e^{-\nu}.$ Therefore
\begin{equation}
|a_3-\frac12a_2^2|^2+|a_2|^2-5\leq\left(2 e^{-2  \nu}  \nu^2+1\right)^2+4 e^{-2  \nu} ( \nu+1)^2-5\equiv G(\nu).\nonumber
\end{equation}
The monotonicity of $G(\nu)$ as a function of $|a_2|$
implies that it holds for
 every  $|a_2|\leq2(\nu+1)e^{-\nu}.$ Indeed
$$ \frac{dG(\nu)}{d\nu}=-16 e^{-4 \nu} \nu^2 \left((\nu-1) \nu+e^{2 \nu}\right)$$
is negative so that  $G(\nu)$ is monotone decreasing and so is $2(\nu+1)e^{-\nu}.$ Hence, $G$ is monotone increasing
in $|a_2|.$ Furthermore, since $G(\nu)$
converges to zero when $\nu$ converges to infinity we conclude that $G$ is a negative function of $|a_2|.$

\vspace{.3cm}
We proceed to find an estimate for $G$ in terms of $|a_2|.$ We consider
$$ J(\nu)=G(\nu)^2-2(2-|a_2|)^3;\;\mbox{where}\; |a_2|=2(\nu+1)e^{-\nu}$$
and show that it is positive for all values of $0\leq|a_2|\leq 2$ or $0\leq\nu\leq \infty.$ We have
\begin{eqnarray}
\frac {J(\nu)}{16e^{-8\nu}}&=&\nu^8-e^{6 \nu} \left(7 \nu^2+10 \nu+5\right)+2 e^{2 \nu} \left(2 \nu^2+2
   \nu+1\right) \nu^4\nonumber\\
   &+&e^{4 \nu} \left(2 \nu^4+8 \nu^3+8 \nu^2+4
   \nu+1\right)+e^{5 \nu} (\nu+1)^3+3 e^{7 \nu} (\nu+1).\nonumber
\end{eqnarray}
To show that $J(\nu)\geq0,$ we distinguish between three cases:

\vspace{.2cm}
\noindent
{\bf Case 1.} $\nu\geq 1.65 \;.$ Note that
$$\frac{J}{e^{6\nu}}\geq - \left(7\nu^2+10 \nu+5\right)+3 e^{\nu}(\nu+1).$$
The proof of the nonegativity of the last expression reduces to showing that
$$K(\nu)=\log[3 e^\nu (\nu+1)/(7\nu^2+10 \nu+5)]\geq0.$$
 Since obviously
$$K'(\nu)=(7 \nu^3+10 \nu^2+\nu)/\left((\nu+1)(7 \nu^2+10 \nu+5)\right)\geq 0$$
and $K(1.65)>0,$ we are done.

\vspace{.3cm}
\noindent
{\bf Case 2.} $0\leq \nu\leq 1 \;.$  We prove that
\begin{equation}H(\nu)=G+\sqrt 2 (2-|a_2|^2)^{\frac32}\leq0\nonumber\end{equation}
for $0\leq\nu\leq1.$ Indeed, $H(0)=0,$ and
\begin{equation}
 H'(\nu)=G'(\nu)-\sqrt2\frac32(2-|a_2|)^{\frac12}(-|a_2|'(\nu)).\nonumber
\end{equation}
To show that $H'(\nu)\leq 0,$  bearing in mind that $G'(\nu)\leq 0,$ we propose to check if
$$G'(\nu)^2-\frac{18}4(2-|a_2|)|a_2|'^2(\nu)\geq0$$
or equivalently
\begin{equation}
4 e^{-8 \nu} \nu^2 \left(64 (\nu-1)^2 \nu^4+128 e^{2 \nu} (\nu-1) \nu^3+64 e^{4 \nu} u^2-9 e^{6 \nu}+9 e^{5 \nu} (\nu+1)\right)\geq0.\nonumber
\end{equation}
Indeed, this is the case, if
$$128 e^{2 \nu} (\nu-1) \nu^3+64 e^{4 \nu} \nu^2-9 e^{6 \nu}+9 e^{5 \nu} (\nu+1)\geq0$$
or dividing by $9e^{2\nu},$ if
\begin{equation}
\frac{128}{9} (\nu-1) \nu^3+\frac{64}{9} e^{2 \nu} \nu^2-e^{3 \nu}(e^\nu-\nu-1)\geq0.
\end{equation}
Invoking the inequality
$$\frac{ e^{3 \nu}(e^\nu-\nu-1)}{\nu^2}=e^{2\nu}e^\nu\frac{e^\nu-\nu-1}{\nu^2}=e^{2\nu}e^\nu\bigg(\frac1{2!}+\frac \nu{3!}+\frac {\nu^2}{4!}+\dots\bigg)\leq e^{2\nu}e(e-1),$$
 it suffices to show that the function
$ \psi_1(\nu)=\frac{e^{2\nu}}{\nu(1-\nu)}\geq 5.83.$ However, the absolute minimum of $\psi_1(\nu)$ in $ [0,1]$ is attained at $\nu=(2-\sqrt{2})/2$ and is larger than 8.

\vspace{.3cm}
\noindent
{\bf Case 3.} $1\leq \nu\leq 1.65 \;.$ Starting from (11) and noting that the first term is nonegative,
it remains to show that
 $$\frac{64}{9} e^{2 \nu} \nu^2-e^{3 \nu}(e^\nu-\nu-1)\geq0.$$
But $\psi_2(\nu)=\frac{e^\nu-\nu-1}{\nu^2}$ is monotone increasing by virtue of its power series expansion
at the origin and hence is bounded in $[1,1.65]$ by 1. Hence it suffices to show that
$e^\nu\leq 64/9$  which is obvious.
This concludes the proof of Lemma 2.

\vspace{.5cm}
\noindent
\section {Proof of Main Theorem }
Let $$f(z)=z+\sum_{n=2}^\infty a_nz^n\in \mathbf{S}$$ and consider the corresponding odd function $f_2(z)\in \mathbf{S}$
$$ f_2(z)=\sqrt{f(z^2)}=\sum_{k=0}^\infty b_k z^{2k+1};\; b_0=1,$$ then
\begin{equation} |a_n|= | \sum_{k=0}^{n-1} b_k b_{n-k}|\leq \sum_{k=0}^{n-1} |b_k|^2.\end{equation}
Furthermore,  from the representation (6) we conclude that
$$\log\left (\sum_{k=0}^\infty b_k z^{2k}\right )= \log\frac {f_2(z)}z =  \log\sqrt {f(z^2)/z^2}=\sum_{k=1}^\infty \gamma_k z^{2k}\equiv h(z^2),$$
so that
$$\exp\big( h(z)  \big)=\sum_{k=0}^\infty b_k z^{k}.$$
We now  apply Lemma 1 in the form $Q_n\leq Q_2$ to the last equality and use Theorem 1 and Remark 1 to get
\begin{eqnarray}
\sum_{k=0}^{n-1} |b_k|^2&\leq& n \frac{(1+|b_1|^2)}2\exp\bigg[M_n(h)-M_2(h)\bigg]\nonumber\\
&\leq& n \frac{1+|b_1|^2}2\exp\bigg[-K_n(h)+\bigg(\frac{1-|\gamma_1|^2}2\bigg)\bigg].\nonumber
\end{eqnarray}
Accordingly, making use of $b_1=-\gamma_1=\frac{a_2}2$ and applying (12), we conclude that
\begin{equation}
 |a_n|\leq \sum_{k=0}^{n-1} |b_k|^2
 \leq n \frac{1+\big|\frac{a_2}2\big|^2}2\exp\bigg[-K_n(h)+\bigg(\frac{1-\big|\frac{a_2}2\big|^2}2\bigg)\bigg]
\end{equation}
Next we differentiate between even and odd $n.$

\vspace{.3cm}
{\bf Case 1: even n}

\vspace{.3cm}
\noindent
Set $n=2m, \; m \geq 2$ in (13) so that
$$|a_{2m}|\leq  2m \frac{1+\big|\frac{a_2}2\big|^2}2\exp\bigg[-K_{2m}(h)+\frac{1-\big|\frac{a_2}2\big|^2}2\bigg].$$
Set $y=\frac{1-\big|\frac{a_2}2\big|^2}2$ so that $1-y=\frac{1+\big|\frac{a_2}2\big|^2}2 .$ In view of
\begin{equation}
(1-x)e^x\leq1-\frac{x^2}2;\;\;\mbox{and}\;\; e^{-x}\leq1-x+x^2/2;\;\; \mbox{for}\;\; x\geq0,
\end{equation}
and invoking Lemma 3 we remain with
\begin{equation}
 |a_{2m}|\leq  2m (1-\frac{y^2}2)\bigg(1-\frac{3}{4(4m^2-1)}y+\bigg(\frac{3}{4(4m^2-1)}\bigg)^2\frac{y^2}2\bigg).\nonumber
 \end{equation}
Ultimately $\frac3{4(4m^2-1)}<1/20$ and $0\leq y\leq1/2$ so that $|a_{2m}|<2m(1-\frac{3y}{4(4m^2-1)})$
and finally
$$2m-|a_{2m}|\geq \frac{6my}{4(4m^2-1)}=\frac{3m}{2(4m^2-1)}\bigg(\frac{4-|a_2|^2}{8}\bigg)\geq\frac{3m}{8(4m^2-1)}(2-|a_2|).$$

\vspace{.3cm}
{\bf Case 2: odd n}

\vspace{.3cm}
\noindent
We set $n=2m+1;\; m\geq 2.$  In accordance with (13) and borrowing $y$ from the previous case we have
$$|a_{2m+1}|\leq (2m +1) (1-\frac{y^2}2)\exp[-K_{2m+1}+y].$$
We start by noting that $1-y^2/2\leq 1-(2-|a_2|)^2/8\equiv 1-t^2/2.$  Then invoking (14) and Lemma 3 we get
\begin{eqnarray}
|a_{2m+1}|&\leq & (2m +1)(1-\frac{ t^2}8)(1-\frac{\sqrt{2}}4 \;K_{2m+1}\;t^{\frac32}+\frac{(\frac{\sqrt2}4 \;K_{2m+1})^2}2\;t^3)\nonumber\\
&\leq& (2m+1)\big(1-\frac{\sqrt2}4 \;K_{2m+1}\;t^{\frac32}-\frac{ t^2}8+\frac{\sqrt2}{32} \;K_{2m+1}\;t^{\frac72}+\frac{1}{16} \;K_{2m+1}^2\;t^3\big).\nonumber
\end{eqnarray}
 Taking into account that $0\leq t\leq2$ we conclude that
\begin{eqnarray}
|a_{2m+1}|&\leq& (2m+1)(1-\frac{\sqrt{2}}4 \;K_{2m+1}\;t^{\frac32}-\frac{ t^2}8(1-K_{2m+1}-K_{2m+1}^2))\nonumber\\
&\leq&(2m+1)(1-\frac{.8}4\; K_{2m+1}\;t^{\frac32}\nonumber\\
       &\;&\;\;\;\;\;\;\;\;\;\;\;\;\;\;\;-\frac{ t^2}8(1-K_{2m+1}-K_{2m+1}^2)-\frac{\sqrt{2}-.8}4\;K_{2m+1} t^{\frac32}).\nonumber
\end{eqnarray}
Noting that
$$-\frac{ t^2}8(1-K_{2m+1}-K_{2m+1}^2))-\frac{\sqrt{2}-.8}4\;K_{2m+1} t^{\frac32}\leq -0.035$$
for $ K_{2m+1}\leq .9,$ which is the case for all $m\geq2,$ we remain with
$$|a_{2m+1}|\leq (2m+1)(1-\frac15 \;K_{2m+1}\;t^{\frac32}).$$
Finally we have that
$$\frac{2m+1-|a_{2m+1}|}{(2-|a_2|)^{\frac32}}\geq \frac{3m(m+1)}{(2m-1)(2m+3)}.$$
This concludes the proof of the main Theorem.

\begin{rk}
 The odd case for $n=3$ can be taken care of by reducing the $.8$ choice in the theorem. A better bound in this case can be
achieved directly from de Branges System as instructed in (\cite{lij} p.167) where this is being taken care of for n=4.
\end{rk}

\vspace{.2cm}
Finally we note that for even $n$ the linear dependence on  $2-|a_2|$ in the main theorem is exact. This is seen via the example
$f_c(z)=\frac z{1-2cz+z^2}.$ For the odd case it is believed to be true.

\vspace {.5cm}

D. Bshouty

daoud@technion.ac.il

D. Aharonov

dova@technion.ac.il

Department of mathematics, Technion IIT, Haifa 32000, Israel.


\begin{thebibliography}{99}

\bibitem{A}  D. Aharonov, Special topics in univalent functions, Lecture Notes, University of Maryland, 1971.

\bibitem{ahd}D.  Aharonov, {\it  On the Bieberbach conjecture for functions with a small second coefficient}, Israel J. Math., 15 (1973), 137–-139.

\bibitem{bil} L. Bieberbach, {\it Über die Koeffizienten derjenigen Potenzreihen, welche eine schlichte Abbildung des Einheitskreises vermitteln},
 S.-B. Preuss. Akad. Wiss., 138 (1916), pp. 940--955.

\bibitem{boe} E. Bombieri,{\it  On the local maximum property of the Koebe function}, Invent. Math., 4 (1967),pp. 26–-67.

\bibitem{dbl} L. de Branges, {\it A proof of the Bieberbach conjecture},Acta Math., 154 (1985),pp. 137–-152.

\bibitem{bsd1} D.  Bshouty, {\it The Bieberbach conjecture for restricted initial coefficients} Math. Z. 182 (1983), pp. 149–-158.

\bibitem{bsd2} D. Bshouty, {\it  A coefficient problem of Bombieri concerning univalent functions}, Proc. Amer. Math. Soc. 91 (1984), pp. 383–-388.

\bibitem{BH} D. Bshouty and W. Hengartner, {\it A variation of the Koebe mapping in a dense subset of S}, Canad. J. Math. 39 (1987), 54–-73.

\bibitem{CS} Z. Charzy\'{n}ski and M. Schiffer,  {\it  A geometric proof of the Bieberbach conjecture for the fourth coefficient}, Scripta Math. 25 (1960), pp. 173–-181.

\bibitem{dxh} X. H. Dong, {\it  A remark on de Branges theorem}, Acta Sci. Natur. Univ. Norm. Hunan., 14 (1991),pp. 193–-197, 205.

\bibitem{DS} P. L. Duren and M. Schiffer, {\it The theory of the second variation in extremum problems for univalent functions}, J. Analyse Math. 10 (1962/1963), pp.  193–-252.

\bibitem{ehg1} G. Ehrig, {\it The Bieberbach conjecture for univalent functions with restricted second coefficients}, J. London Math. Soc. (2) 8 (1974),pp. 355–-360.
\bibitem{fic1} C. H.  FitzGerald, {\it Quadratic inequalities and coefficient estimates for schlicht functions}, Arch. Rational Mech. Anal. 46 (1972), pp. 356–-368.

\bibitem{fic2} C. Fitzgerald,  {\it The Bieberbach conjecture: retrospective}, Notices Amer. Math. Soc. 32 (1985), pp. 2–-6.

\bibitem{GS1} P. R. Garabedian and M. Schiffer,   {\it  A proof of the Bieberbach conjecture for the fourth coefficient}, J. Rational Mech. Anal. 4 (1955),  pp. 427–-465.

\bibitem{GS} P. R. Garabedian and M. Schiffer, {\it The local maximum theorem for the coefficients of univalent functions}, Arch. Rational Mech. Anal. 26 (1967), pp. 1–-32.

\bibitem{GRS} P. R. Garabedian, G. G. Ross and M. M. Schiffer, {\it  On the Bieberbach conjecture for even n}, J. Math. Mech. 14 (1965) 975–-989.

\bibitem{gol}  G. M.  Goluzin, {\it Geometric theory of functions of a complex variable}, Translations of Mathematical Monographs, Vol. 26, American Mathematical Society, Providence, R.I. (1969).

\bibitem{gor} J. G\'{o}rski, {\it A certain minimum problem in the class S}, Ann. Univ. Mariae Curie-Skłodowska, Sect. A 22-24(1968-1970),  73--77.

\bibitem{grr} R. Greiner and O. Roth, {\it On support points of univalent functions and a disproof of a conjecture of Bombieri}, Proc. Amer. Math. Soc. 129 (2001), pp. 3657–-3664.

\bibitem{grh} H.  Grunsky,  {\it  Koeffizientenbedingungen f\"{u}r schlicht abbildende meromorphe Funktionen},
     Math. Z. 45 (1939), pp. 29–-61.

\bibitem{had} D. Hamilton, Oral Communication.

\bibitem{haw}  W. K. Hayman, {\it Research problems in function theory}, The Athlone Press, University of London, London (1967).

\bibitem{huj} J. A. Hummel, {\it The coefficients of starlike functions}, Proc. Amer. Math. Soc. 22 (1969), pp. 311–-315.

\bibitem{lij} J.-L. Li, {\it Notes on the Duren-Leung conjecture}, J. Math. Anal. Appl. 332 (2007), pp. 164--170.

\bibitem{lok} K. L\"{o}wner, {\it  Untersuchungen \"{u}ber schlichte konforme Abbildungen des Einheitskreises I},
Math. Ann. 89 (1923), pp. 103–-121.

\bibitem{mii}  I. M. Milin,{\it The coefficients of schlicht functions}, (Russian) Dokl. Akad. Nauk SSSR 176 (1967),pp. 1015–-1018. English transl. [Soviet Math. Dokl. 8 (1967), pp. 1255–-1258].

\bibitem{M}  I. M. Milin, Univalent functions and orthonormal systems,   Izdat. ``Nauka'', Moscow, 1971. English transl. [ Math. Monos., vol. 50, Amer. Math. Soc., Providence, R.I., 1977.]

\bibitem{P}  Ch. Pommerenke, Univalent functions, Vandenhoeck and Ruprecht, Göttingen, 1975.

\bibitem{PR} D. Prokhorov and O. Roth, {\it On the local extremum property of the Koebe function}, Math. Proc. Cambridge Philos. Soc., 136 (2004), pp. 301--312.

\bibitem{PV} D. Prokhorov and A. Vasilʹev, {\it Optimal control in Bombieri's and Tammi's conjectures}, Georgian Math. J., 12 (2005), pp. 743--761.

\bibitem{scm2}M. Schiffer, {\it  Sur un probl\`{e}me d'extr\'{e}mum de la repr\'{e}sentation conforme},  Bull. Soc. Math. France 66 (1938),pp. 48–-55.

\bibitem{scm1} M. Schiffer, {\it A Method of variation within the family of simple functions}, Proc. London Math. Soc. S2-44 (1938), pp.  432--449.

\end{thebibliography}
\end{document}